\documentclass[reqno,11pt]{amsart}

\newtheorem{rem}{Remark}

\newtheorem{prop}{Proposition}
\newtheorem{cor}{Corollary}
\newtheorem{df}{Definition}

\newcommand\eps\varepsilon
\newcommand\ph\varphi
\newcommand\kap\varkappa

\usepackage{enumerate}

\begin{document}

\title[Evolution Differential Equations in Fr\'echet Space]
{ Evolution Differential Equations in Fr\'echet Space with Schauder Basis}

\author[Oleg Zubelevich]{Oleg Zubelevich\\ \\\tt
 Dept. of Theoretical mechanics,  \\
Mechanics and Mathematics Faculty,\\
M. V. Lomonosov moscow State University\\
Russia, 119899, moscow, Vorob'evy gory, MGU \\
 }
\date{}
\thanks{Partially supported by grants
 RFBR  12-01-00441.,  Science Sch.-2964.2014.1}

\begin{abstract}We consider evolution differential equations in Fr\'echet spaces that possess unconditional Schauder basis and construct a version of the  majorant functions method   to obtain existence theorems for Cauchy problems.  Applications to PDE and ODE have been considered.
\end{abstract}

\maketitle\tableofcontents
\numberwithin{equation}{section}
\newtheorem{theorem}{Theorem}[section]
\newtheorem{lemma}[theorem]{Lemma}
\newtheorem{definition}{Definition}[section]

\section{Introduction}Countable systems of ordinary differential equations appear in different areas of differential equations and applications, see, for example \cite{sam},\cite{kuk}. 

The most famous  problem which leads to such an object is the Cauchy-Kovalevskaya problem in nonanalytic in time statement. To reduce this problem to the countable system of ODE one must expand the solution to the Taylor series in spacial variables and   substitute this expansion to the corresponding initial value problem then  the Taylor coefficients satisfy infinite system of ODE. 

The  Cauchy-Weierstrass-Kovalevskaya method of majorant functions can be modified for nonanalytic in time statement to obtain corresponding existence theorem \cite{zu}. Generally, being applied to Cauchy-Kovalevskaya problem, this modification does not give anything different from the  results of Nirenberg and Nishida  \cite{nir}.
Nevertheless,  in some cases this method allows to obtain global in time existence theorems or at least effective estimates for the solution's existence time \cite{tr_zu}. 

Another application of the majorant functions method is the initial value problems with non-Lipschitz right hand side. It is well known that in infinite dimensional space such problems in general do not have solutions. But the majorant functions method allows to prove the existence theorems in some special cases. 

This article is devoted to the generalisation of this method for countable systems of ODE in the Fr\'echet spaces that possess the Schauder basis. 

For example,   $\mathcal{D}(\mathbb{T}^m),\quad \mathbb{T}^m=\mathbb{R}^m/(2\pi\mathbb{Z})^m$ is a Fr\'echet space with the unconditional Schauder basis $\{e^{i(k,x)}\},\quad k\in\mathbb{Z}^m.$ Other examples see below.
\section{ Main Theorems}
Let $E$ stand for a  Fr\'echet space. Its topology is defined by the collection of seminormes $\{\|\cdot\|_n\}_{ n\in\mathbb{N}}$.

 Recall that such a space is completely metrizable by the following metrics 
$$\rho(x,y)=\sum_{k=1}^\infty \frac{1}{2^k}\min\{1,\|x-y\|_k\}.$$ 

\begin{df}A sequence $\{e_k\}_{k\in\mathbb{N}}\subset E$ is called a Schauder basis in $E$ if for every $x\in E$ there is a unique sequence of  scalars $\{x_k\}_{k\in\mathbb{N}}$ so that
\begin{equation}\label{cdvfghyu8765fghjknhbgv}x=\sum_{k=1}^\infty x_ke_k.\end{equation}
This series is convergent in the topology of $E$.

We shall say that  $\{e_k\}_{k\in\mathbb{N}}$ is an unconditional basis if for any $x\in E$ and for any permutation $\pi:\mathbb{N}\to\mathbb{N}$ the sum
$$ \sum_{k=1}^\infty x_{\pi(k)}e_{\pi(k)}$$ is convergent.\end{df}
In the sequel we  assume that $E$  possesses an unconditional Schauder basis.

Introduce a notation $I_T=[0,T],\quad T>0$. By definition for $T=\infty$ put $I_\infty=[0,\infty)$.
If it is not explicitly  specified that $T=\infty$ , we assume that $T<\infty$.

\begin{df}We shall say that an element $y=\sum_{k=1}^\infty y_ke_k$ is a majorant for an element $x=\sum_{k=1}^\infty x_ke_k$ and write $x\ll y$ iff$$ |x_k|\le y_k,\quad k\in\mathbb{N}.$$\end{df}

\begin{df}We shall say that $x(t)\in C^1(I_T,E)$ iff for each $t\in I_T$ there exists an element $\dot x(t)$ such that for all $i$ one has
\begin{equation}\label{sfsfdddddddddddd}\lim_{h\to 0}\Big\|\frac{x(t+h)-x(t)}{h}-\dot x(t)\Big\|_i=0.\end{equation}
And the element $\dot x$ belongs to $C(I_T,E)$.

In formula (\ref{sfsfdddddddddddd}) it is assumed that if $t=0$ then $h>0$ and $h<0$ provided $t=T$.\end{df}

Fix an element $y\in E$ and let $\mathcal{ X}_j[y]:E\to E$ stand for the following affine mappings $\mathcal X_1[y]x= y_1e_1+\sum_{k=2}^\infty  x_ke_k,$ $$\mathcal X_j[y]x= \sum_{k=1}^{j-1}x_ke_k+y_je_j+\sum_{k=j+1}^\infty  x_ke_k,\quad j>1.$$

Let a function $X(t)=\sum_{k=1}^\infty X_k(t)e_k\in C(I_T,E)$ be such that $$X_k(t)\ge 0,\quad k\in \mathbb{N},\quad t\in I_T,$$
and $X_k(t)\in C^1(I_T)$.

Introduce a set $$W_X=\{(t,x)\in I_T\times E\mid x\ll X(t)\}.$$

Consider the following initial value problem
\begin{align}
\dot x&=f(t,x),\quad x(0)=\hat x,\label{qq23}\\  f(t,x)&=\sum_{k=1}^\infty  f_k(t,x)e_k,\quad f\in C(W_X,E).\nonumber \end{align}

\begin{theorem}\label{re356qqbyh6}
Suppose that $X_k(t)>0,\quad k\in\mathbb{N},\quad t\in I_T$ and 
 for each $(t,x)\in W_X$ one has
$$\pm f_k(t,\mathcal{X}_k[\pm X(t)]x)\le \dot X_k(t),\quad\hat x\ll X(0).$$
(Here and in the sequel this means that for each $k$ two inequalities hold.)

Then problem (\ref{qq23}) has a solution  $x(t)\in C^1(I_T,E)$ such that $$x(t)\ll X(t),\quad t\in I_T.$$ \end{theorem}\begin{rem}The function $X(t)$ that satisfies the conditions of Theorem \ref{re356qqbyh6} is called a  majorant function for  problem (\ref{qq23}).\end{rem} This theorem develops corresponding results of \cite{zu} and, like the theorems from that article, implies the classical Cauchy-Kovalevskaya theorem and a number of its generalisations.

\begin{theorem}\label{fdfgdfgtt} Suppose that $T=\infty$ and the function $f$ is $\omega-$periodic ($\omega>0$) in $t$. 

Suppose also that $X_k(t)>0,\quad k\in\mathbb{N},\quad t\in I_T$ and 
 for each $(t,x)\in W_X$ one has
$$\pm f_k(t,\mathcal{X}_k[\pm X(t)]x)\le \dot X_k(t),$$
and $X(\omega)\ll X(0).$

Then problem (\ref{qq23}) has a solution  $\tilde x(t)\in C^1(I_\infty,E)$ such that $$\tilde x(t)\ll X(t),\quad \tilde x(t+\omega)=\tilde x(t)\quad t\in I_\infty.$$ \end{theorem}Theorems \ref{re356qqbyh6} and \ref{fdfgdfgtt} are proved in Section \ref{fdgv}.

The following technical proposition is useful for proving continuity of some mappings. 
\begin{prop}\label{fg555} Let $A=\sum_{k=1}^\infty A_ke_k\in E,\quad A_k\ge 0$ be a fixed element.

Assume that a sequence $x_n=\sum_{k=1}^\infty x_{kn}e_k$ belongs to $$K_A=\Big\{y=\sum_{k=1}^\infty y_ke_k\in E\mid | y_k|\le A_k\Big\}$$ and this sequence is weakly convergent: for all $k$ it follows that $x_{kn}\to x_k$ as $n\to\infty$. Then $x=\sum_{k=1}^\infty x_{k}e_k\in K_A$ and the sequence  is  convergent in $E$ i.e. $\rho(x_n,x)\to 0.$\end{prop} It is proved by the methods developed in Section \ref{fdgv}.

\subsection{Non-negative Solutions} In this section we formulate another pair of theorems which  belong to the same range of ideas. We do not bring their proofs since they repeat the argument of Section \ref{fdgv} up to evident modifications.

Endow the space $E$ with partial order $\prec$ by the following rule.
\begin{df}We shall  write $x=\sum_{k=1}^\infty x_ke_k\prec y=\sum_{k=1}^\infty y_ke_k$ iff$$ x_k\le y_k,\quad k\in\mathbb{N}.$$\end{df}
Introduce a set $$W_X^+=\{(t,x)\in I_T\times E\mid 0\prec x\prec X(t)\}.$$
 
Assume  that $f\in C(W_X^+,E)$.
\begin{theorem}\label{re31156qqbyh6}
Suppose that $X_k(t)>0,\quad k\in\mathbb{N},\quad t\in I_T$ and 
 for each $(t,x)\in W_X^+$ one has
$$  f_k(t,\mathcal{X}_k[ X(t)]x)\le \dot X_k(t),\quad 0\prec\hat x\prec X(0)$$
and $f_k(t,\mathcal{X}_k[0]x)\ge 0.$

Then problem (\ref{qq23}) has a solution  $x(t)\in C^1(I_T,E)$ such that $$0\prec x(t)\prec X(t),\quad t\in I_T.$$ \end{theorem}

\begin{theorem}\label{fdfg1dfgtt} Suppose that $T=\infty$ and the function $f$ is $\omega-$periodic ($\omega>0$) in $t$. 

Suppose also that $X_k(t)>0,\quad k\in\mathbb{N},\quad t\in I_T$ and 
 for each $(t,x)\in W_X^+$ one has
$$f_k(t,\mathcal{X}_k[ X(t)]x)\le \dot X_k(t),$$
and  $f_k(t,\mathcal{X}_k[0]x)\ge 0.$ Moreover suppose that
 $X(\omega)\ll X(0).$

Then problem (\ref{qq23}) has a solution  $\tilde x(t)\in C^1(I_\infty,E)$ such that $$0\prec\tilde x(t)\prec X(t),\quad \tilde x(t+\omega)=\tilde x(t)\quad t\in I_\infty.$$ \end{theorem}

\section{Applications}

\subsection{ Linear PDE}To release our exposition from technical details we restrict ourselves to the case of PDE with one-dimensional spatial variable. However  considered below propositions can easily be obtained for corresponding systems with multidimensional spatial variable.

\subsubsection{The Existence Theorem}

Let $\mathcal O(\mathbb{C})$ stand for the space of entire functions $u:\mathbb{C}\to\mathbb{C}$. This is a  Fr\'echet space with seminorms $$\|v\|_n=\max_{|z|\le n}|v(z)|,\quad n\in\mathbb{N}$$ and the Schauder basis is $e_j=z^j,\quad j\in\mathbb{Z}_{+}=\{0,1,2,\ldots\}.$

For the space $E\subset\mathcal O(\mathbb{C})$  take the space of entire functions $$v: \mathbb{C}\to \mathbb{C}$$ such that $ \overline{v(z)}=v(\overline z).$

Fix arbitrary positive number $T$  and take  functions $a(t),b(t)\in C(I_T,\mathbb{R})$. Consider the following initial value problem
\begin{equation}\label{5vb54}v_t(t,z)=b(t)v(t,z)+a(t)z^m\frac{\partial^N v(t,z)}{\partial z^N},\quad v(0,z)=\hat v(z).\end{equation}
Here $N>m\ge 0$ are some  integers.

Introduce the following notations $$q_{jmN}=\frac{(j-m+N)!}{(j-m)!},\quad j\ge m,\quad a^*=\|a\|_{C(I_T)}.$$ We assume that $a^*\ne 0$.
Then  take arbitrary positive constants $U_0,\ldots U_{N-1}$ and define other constants recurrently   
$$U_{j-m+N}=\frac{U_j}{a^*q_{jmN}},\quad j\ge m.$$
It is not hard to show that $U(z)=\sum_{j=0}^\infty U_jz^j\in E$.

\begin{prop}\label{dfgdfgd44}Suppose that $\hat v\ll U$. Then problem (\ref{5vb54}) has a solution $v(t,z)\in C^1(I_T,E)$ and $v(t,z)\ll e^{\int_0^tb(s)ds+ t}U(z)$ for all $t\in I_T$.\end{prop}
Note that this proposition does not follow from results of \cite{Dubinskii}.

Indeed, after the change of function $v=e^{\int_0^tb(s)ds+ t}u$
problem (\ref{5vb54}) takes the form
\begin{equation}\label{rgery5dbb}
u_t(t,z)=-u(t,z)+a(t)z^m\frac{\partial^N u(t,z)}{\partial z^N},\quad u(0,z)=\hat v(z).
\end{equation}
In coordinate notation problem (\ref{rgery5dbb}) has the form $u(t,z)=\sum_{j=0}^\infty u_j(t)z^j$,$$
\dot u_j=-u_j,\quad u_j(0)=\hat v_j,\quad j<m,$$
and
\begin{equation}\label{dbb}\dot u_j=-u_j+q_{jmN}a(t)u_{j-m+N},\quad u_j(0)=\hat v_j,\quad j\ge m.\end{equation}
To apply Theorem \ref{re356qqbyh6} to problem (\ref{dbb}) observe that
for any $$|u_l|\le U_l,\quad l\ge m$$ we have
\begin{align}\pm\big(-(\pm U_j)&+a(t)q_{jmN}u_{j-m+N}\big)\nonumber\\
&\le
-U_j+a^*q_{jmN}U_{j-m+N}=0=\dot U_j,\quad j\ge m.\nonumber\end{align}

The Proposition is proved.

\subsubsection{Periodic Solutions}
Let us redefine sequence $\{U_k\}$.
Take a sequence $F_k\ge 0,\quad k\in\mathbb{N}$ and let  $U_0,\ldots U_{N-1}$ be arbitrary positive constants. Then put
$$U_{j-m+N}=\frac{U_j}{a^*q_{jmN}+F_j},\quad j\ge m.$$
It is not hard to show that $U(z)=\sum_{j=0}^\infty U_jz^j\in E$.

Consider the following system
\begin{equation}\label{ddrege555}u_t(t,z)=-u(t,z)+a(t)z^m\frac{\partial^N u(t,z)}{\partial z^N}+f(t,z).
\end{equation}
Assume that the function
$$f(t,z)=\sum_{k=0}^\infty f_k(t)z^k\in C(I_\infty,E)$$ and $f_k$ are the $\omega-$periodic functions.

\begin{prop}\label{ewffwef}Suppose that for $j\ge m$ it follows that
$$\max_{t\in I_\omega}|f_j(t)|\le F_j U_{j-m+N}.$$
Then system (\ref{ddrege555}) has an $\omega-$periodic solution $u(t,z)
\in C^1(I_\infty,E).$\end{prop}

In coordinate notation problem (\ref{ddrege555}) has the form $$
\dot u_j=-u_j+f_j,\quad j<m,$$
and
\begin{equation}\label{dasfbb}\dot u_j=-u_j+q_{jmN}a(t)u_{j-m+N}+f_j,\quad  j\ge m.\end{equation}
To apply Theorem \ref{fdfgdfgtt} to problem (\ref{dasfbb}) observe that
for any $$|u_l|\le U_l,\quad l\ge m$$ we have
\begin{align}\pm\big(-(\pm U_j)&+a(t)q_{jmN}u_{j-m+N}+f_j(t)\big)\nonumber\\
&\le
-U_j+a^*q_{jmN}U_{j-m+N}+F_jU_{j-m+N}=0=\dot U_j,\quad j\ge m.\nonumber\end{align}

The Proposition is proved.

\subsection{Periodic Solutions to the Smoluchowski Equation }

In this section we consider the following IVP
\begin{equation}\label{sawq}
\dot x_k=c_k+\frac{1}{2}\sum_{i+j=k}b_{ij}x_ix_j-x_k\sum_{j}b_{kj}x_j,\quad x_k(0)=\hat x_k,\quad i,j,k\in\mathbb{N}.\end{equation}
The functions $c_i(t),b_{ij}(t)\in C(I_T)$ are non negative valued, $$b_{ij}=b_{ji},\quad \hat x_k\ge 0.$$ For this IVP the non negative solutions $x_k(t)\ge 0$ are of interest.

In \cite{McLeod},\cite{wite} the existence theorems have been proved under the following assumptions $b_{ij}(t)\le (i+j)^\alpha,\quad \alpha\in [0,1]$ and $\hat x_k,c_k$ are decreased as $c/k^p$ with some $p>\alpha$. The obtained solutions are bounded in certain norm on   bounded intervals. 

 We do not assume anything about growth of coefficients $b_{ij}$ for $i\ne j$, but our assumption on growth of coefficients $b_{kk}$ is strong enough. 
Under these assumptions we  prove the existence of bounded for all time solutions and the existence  of a periodic solution when the coefficients $b_{ij},c_k$ are periodic. 

Let us put $$c_k(t)\le C_k,\quad b_{kk}(t)\ge\beta_k,\quad b_{ij}(t)\le B_{ij}.$$
In these formulas $i,j,k\in\mathbb{N}$ and the inequalities hold for all $t\in I_T$ with some non negative constants $C_k,B_{ij},\beta_k$. 

Introduce a sequence
$$X_1=\sqrt{C_1/\beta_1},\quad X_k=\sqrt{\frac{C_k+\frac{1}{2}\sum_{i+j=k}B_{ij}X_iX_j}{\beta_k}},\quad k=2,3,\ldots.$$ We assume that the constants $\beta_k$ are large such that

 $$b_k=\sum_{j=1}^\infty B_{kj}X_j<\infty.$$

Let us put
$$A_k=C_k+\frac{1}{2}\sum_{i+j=k}B_{ij}X_iX_j+X_kb_k,\quad F_k=\max\{1,A_k,X_k\}.$$

Introduce a Banach space $E$ of sequences $x=\{x_k\}$ with the following norm
$$\|x\|=\sum_{k=1}^\infty \frac{1}{k^2F_k}|x_k|<\infty.$$ Evidently, this space possesses an unconditional Schauder basis $e_j=\{\delta_{ij}\}_{i\in\mathbb{N}}.$
Note that $$A=\sum_{k\in\mathbb{N}}A_ke_k,\quad X=\sum_{k\in\mathbb{N}}X_ke_k\in E.$$

\begin{prop}\label{rrerf}
Assume that $\hat x\prec X$.

Then problem (\ref{sawq}) has a solution $x(t)\in C^1(I_T,E)$ such that $$0\le x_k(t)\le X_k,\quad t\in I_T.$$ 

Moreover, if  the  functions $b_{ij},c_j$ are $\omega-$ periodic then there is a solution $\tilde x(t)\in C^1(I_\infty,E),\quad \tilde x_k(t)\ge 0$ that is also $\omega-$periodic and $0\le \tilde x_k(t)\le X_k$. \end{prop}

\proof So we wish to apply  theorems \ref{re31156qqbyh6}, \ref{fdfg1dfgtt}.

For $0\le x_s\le X_s,\quad s\in\mathbb{N}$ and $x_k=0$ the condition of the theorems is satisfied identically
$$
 c_k+\frac{1}{2}\sum_{i+j=k}b_{ij} x_ix_j\ge 0.$$

Another condition of the theorems to check is for  $0\le x_s\le X_s,\quad s\in\mathbb{N}$ and $x_k=X_k:$

\begin{align}c_k&+\frac{1}{2}\sum_{i+j=k}b_{ij}x_ix_j-x_k\sum_{j}b_{kj}x_j\nonumber\\
&\le C_k+\frac{1}{2}\sum_{i+j=k}B_{ij}X_iX_j-X_k^2 \beta_k\le \dot X_k=0\nonumber\end{align}

It remains to show that the mapping
$$f(t,x)=\sum_{k=1}^\infty f_k(t,x)e_k,\quad f_k=c_k+\frac{1}{2}\sum_{i+j=k}b_{ij}x_ix_j-x_k\sum_{j}b_{kj}x_j$$
is continuous as a mapping of $W_X^+$ to $E$. 
Observe that for each $(t,x)\in W_X^+$ we have $f(t,x)\ll A$.
Now the continuity follows from Proposition \ref{fg555}.

\section{Proof of Main Theorems }\label{fdgv}
\subsubsection{A Short Digression in Functional Analysis}

Let $\mathcal P_n:E\to E$ be the projection
$$\mathcal P_n\Big(\sum_{k=1}^\infty x_ke_k\Big)=\sum_{k=1}^n x_ke_k.$$
Let us also put  $\mathcal Q_n=\mathrm{id}-\mathcal P_n.$

\begin{theorem}\label{e5b5j}Let  
$\lambda=\{\lambda_j\}_{j\in\mathbb{N}}\in \ell_\infty$ and let $$\mathcal M_\lambda x=\sum_{k=1}^\infty \lambda_k x_k e_k.$$
Then for any number $i'$ there exist a number i and a positive constant $c$ both independent on $\lambda$ such that 
$$\|\mathcal M_\lambda x\|_{i'}\le c\|\lambda\|_{\infty}\cdot\|x\|_i.$$
\end{theorem}Particularly, Theorem \ref{e5b5j} implies that the operators $\mathcal P_n,\mathcal Q_n$ are continuous. However this fact needs an independent proof since Theorem \ref{e5b5j}  is based  upon it by itself. 

Theorem \ref{e5b5j} and the continuity of the projections are proved in  Section \ref{rvbere6g}.

\begin{lemma}\label{dfgg}The set $W_X$ is a compact set in $I_T\times E$.\end{lemma}
\proof  

Consider continuous mappings $$v_n:I_T\to E,\quad v_n(t)=\mathcal Q_nX(t).$$
This sequence is pointwise convergent to zero: $v_n(t)\to 0,\quad n\to\infty$ for any fixed $t\in I_T$. On the other hand this sequence is  uniformly continuous on  $I_T$.

Indeed, by Theorem \ref{e5b5j} the mappings $\mathcal Q_n$ are uniformly continuous thus for any $i'$ there exist a constant $c>0$ and a number $i$ such that $$\|v_n(t')-v_n(t'')\|_{i'}\le c\|X(t')-X(t'')\|_{i}.$$ But the mapping $X$ is uniformly continuous on the compact set $I_T$.

 Consequently, $v_n(t)\to 0$ uniformly in $I_T$ \cite{Shvarz}.

Evidently,  the set $W_X$ is closed. We prove the Lemma if show that the sets $$A_n=\{(t,\mathcal P_n x)\in I_T\times E\mid x\ll X(t)\}$$ form a sequence of compact $\epsilon-$nets in $W_X$. 

Indeed, each set  $A_n$ is contained in $\mathbb{R}^{n+1}$, closed and bounded.

Let us take an element $(t,x)\in W_X$; and employ  Theorem \ref{e5b5j} with $$\lambda_k(t)=x_k/X_k(t),\quad |\lambda_k(t)|\le 1$$ then it follows that  for any number $i'$ there exist a number $i$ and a constant $c$ such that
$$\|x-\mathcal P_nx\|_{i'}=\|\mathcal M_{\lambda(t)}\mathcal Q_nX(t)\|_{i'}\le c\|v_n(t)\|_i.$$
Therefore $ \sup_{(t,x)\in W_X}\|x-\mathcal P_nx\|_{i'}\to 0$ as $n\to\infty$. 

The Lemma is proved.

By the analogous argument one obtains the following lemma.

\begin{lemma}\label{rb45y}Take an element $U=\sum_{k=1}U_ke_k,\quad U_k\ge 0$. Then $$K_U=\{u\in E\mid u\ll U\}$$ is a compact set.\end{lemma}

\begin{theorem}[Arzela, Ascoli, \cite{Shvarz}]\label{coppp}Consider a   set $K\subset C(I_T,E)$. Suppose that 

1) for any $t\in I_T$ the set $K_t=\{x(t)\mid x(\cdot)\in K\}\subset E$ is compact.

2) for any $\epsilon>0$ and for any $n\in\mathbb{N}$ there exist a constant $\delta>0$ such that if $t',t''\in I_T,\quad |t'-t''|<\delta$ then $$\|x(t')-x(t'')\|_n\le \epsilon.$$  Then $K$ is a compact set.
\end{theorem}

\subsubsection{Back to Proof of Theorem \ref{re356qqbyh6}}

We approximate problem (\ref{qq23}) by the following finite dimensional ones
\begin{equation}\label{ev5b5}\dot y^n=\mathcal P_nf(t,y^n),\quad y^n(0)=\hat y^n=\mathcal P_n \hat x,\quad y^n=\sum_{j=1}^ny_je_j.\end{equation}
By Theorem \ref{re356byh6} all the problems (\ref{ev5b5}) have solutions $y^n(t)\in C^1(I_T,\mathbb{R}^n)$ and \begin{equation}\label{dgbdyt}(t,y^n(t))\in W_X,\quad t\in I_T.\end{equation}

By Theorem \ref{e5b5j} and Lemma \ref{dfgg} for any $i'$ there is a number $i$ and a constant $c$ such that 
\begin{align}\sup\{\|&\dot y^n(t)\|_{i'}\mid n\in\mathbb{N},\quad t\in I_T\} \nonumber\\
&\le \sup\{\|\mathcal P_nf(t,x)\|_{i'}\mid (t,x)\in W_X,\quad n\in\mathbb{N}\}\nonumber\\&\le c\sup_{(t,x)\in W_X}\|f(t,x)\|_{i}
\le C_{i'}<\infty.\nonumber\end{align}
For any $t',t''\in I_T$ this implies
$$\|y^n(t')-y^n(t'')\|_{i'}=\Big\|\int_{t'}^{t''}\dot y^n(s)ds\Big\|_{i'}\le C_{i'}|t'-t''|.$$
By Theorem \ref{coppp} and Lemma \ref{rb45y} the sequence $\{y^n\}$ contains a subsequence that is convergent in $C(I_T,E)$. Denote this subsequence in the same manner: $$y^n(\cdot)\to x(\cdot)\quad \mbox{in} \quad C(I_T,E).$$Since the operators $\mathcal P_n$ are continuous, formula (\ref{dgbdyt}) implies  $x(t)\ll X(t),\quad t\in I_T$.

Our next goal is to show that $x(t)$ is the desired solution to problem  (\ref{qq23}).
Rewrite problem (\ref{ev5b5}) as follows
\begin{equation}\label{rebsfdsdf}
y^n(t)-\hat y^n=\int_0^t\mathcal P_nf(s,y^n(s))ds.\end{equation}
Passing to the limit as $n\to \infty$ in the left side of this formula we obtain $x(t)-\hat x$.

Consider the right hand side of formula (\ref {rebsfdsdf}). 
\begin{lemma}\label{db6b7}For all $i\in\mathbb{N},\quad t\in I_T$ one has$$\Big\|\int_0^t\mathcal P_nf(s,y^n(s))ds-\int_0^tf(s,x(s))ds\Big\|_i\to 0$$ as $n\to\infty.$ 

The integrals are  understood in the sense of  Millionshchikov \cite{mill}.\end{lemma}\proof
Estimate this expression by parts
\begin{align}\Big\|\int_0^t&\mathcal P_nf(s,y^n(s))ds-\int_0^tf(s,x(s))ds\Big\|_i\nonumber\\
&\le \int_0^t \|\mathcal P_n(f(s,y^n(s))-f(s,x(s)))\|_ids\nonumber\\
&+\int_0^t \|\mathcal Q_n f(s,x(s))\|_ids.\nonumber
\end{align}
Then due to Theorem \ref{e5b5j} we have 
$$\|\mathcal P_n(f(s,y^n(s))-f(s,x(s)))\|_i\le c_i\|f(s,y^n(s))-f(s,x(s))\|_{i'}\to 0.$$ Since the function $f$ is uniformly continuous in the compact set $W_X$, this limit is uniform in $s\in I_T$. 

The set $f(W_X)$ is a compact set as an image of a compact set. The operators $\mathcal Q_n$ are uniformly continuous (Theorem \ref{e5b5j}). Consequently, the convergence  $\mathcal \|Q_n f(s,x(s))\|_i\to 0$ is uniform in $s\in I_T$ \cite{Shvarz}. 

The Lemma is proved.

From Lemma \ref {db6b7} and formula (\ref{rebsfdsdf}) it follows that
$$
x(t)-\hat x=\int_0^tf(s,x(s))ds.$$
Consequently $x(t)\in C^1(I_T,E)$ and $\dot x(t)=f(t,x(t)).$ \cite{mill}.

In coordinate notation this implies
$$ x_k(t)-\hat x_k=\int_0^tf_k(s,x(s))ds$$ or
$$\dot x_k(t)=f_k(t,x(t)),\quad x_k(0)=\hat x_k.$$
This particularly implies that the series $\sum_{k=1}^\infty \dot x_k(t)e_k$ is convergent for each $t$.

Theorem \ref{re356qqbyh6} is proved.

\subsubsection{Proof of Theorem \ref{fdfgdfgtt}}

 By Theorem \ref{5evt5sdds} all the problems (\ref{ev5b5}) have $\omega-$periodic solutions $\tilde y^n(t)$ such that 
$$ \tilde y^n(t)\ll X(t).$$ By the same argument as above, the set $\{\tilde y^n(\cdot)\}$ is  relatively compact in $C(I_\omega,E)$. Let $y_*(t)$ be an  accumulation point of this set. Then the function $$\tilde x(t)=y_*(\tau),\quad \tau\in I_\omega,\quad  t=\tau\pmod{\omega}$$ is the periodic solution.

The Theorem is proved.

\section{Finite Dimensional Case}

\subsection{Estimates From Above}
In this section we consider ordinary differential equations in $\mathbb{R}^m.$
Introduce several notations
$$\mathbb{R}^m_+=\{x=(x_1,\ldots, x_m)\in \mathbb{R}^m\mid x_k\ge 0,\quad k=1,\ldots,m\},\quad I_T=[0,T].$$
We shall say that a vector $X=(X_1,\ldots,X_m)\in\mathbb{R}^m_+$ majorizes  a vector $x=(x_1,\ldots,x_m)\in\mathbb{R}^m$   iff
$$|x_k|\le X_k,\quad k=1,\ldots, m.$$ This relation is written as $x\ll X$.
Suppose a function $X(t)\in C^1(I_T,\mathbb{R}^m)$ is such that
 $$X(t)\in \mathbb{R}^m_+,$$
for all $t\in I_T$ with some fixed $T>0$.

Let $U\subset I_T\times\mathbb{R}^m$ be an open neighbourhood    of the  following set
$$W_X=\{(t,x)\in I_T\times\mathbb{R}^{m}\mid x\ll X(t)\}.$$

\subsubsection{Lipschitz Case}\label{5g6rryryrty}

Introduce a  function $f(t,x)\in C(U,\mathbb{R}^m)$ which is a locally Lipschitz function in the second argument. In short words $f$ is a function such that 
the initial value problem 
\begin{equation}\label{bh7}
\dot x=f(t,x),\quad x(0)=\hat x,\quad x=(x_1,\ldots, x_m).\end{equation}
satisfies the standard Cauchy existence and uniqueness theorem in $U$.

\begin{theorem}\label{rebyh6}
Suppose that $X_k(t)>0,\quad k=1,\ldots,m,\quad t\in I_T$ and 
 for each $(t,x)\in W_X$ one has
\begin{equation}\label{e5tb54}\pm f_k(t,x_1,\ldots,x_{k-1},\pm X_k(t),x_{k+1},\ldots,x_m)\le \dot X_k(t),\quad\hat x\ll X(0).\end{equation}

Then problem (\ref{bh7}) has a solution  $x(t)\in C^1(I_T,\mathbb{R}^m)$ such that
 $x(t)\ll X(t)$. \end{theorem}

\subsubsection{Proof of Theorem \ref{rebyh6}}\label{5by45}
Let $\psi(s)\in C^\infty(\mathbb{R}),\quad \mathrm{supp}\,\psi\subset[-1,1]$ be a non-negative function, $\psi(0)=1$. Construct a function $$\phi_\epsilon(s)=\epsilon\psi(s/\epsilon),\quad \epsilon>0.$$
Choose $\epsilon_0>0$  such that the equalities
\begin{equation}\label{e5g6}\phi_\epsilon(\pm 2X_k(t))=0,\quad k=1,\ldots,m\end{equation}hold
 for all $\epsilon\in(0,\epsilon_0),\quad t\in I_T$.

Define a function  $f^\epsilon$  as follows
  $$f^\epsilon_k(t,y)= f_k(t,y)-\phi_\epsilon(y_k-X_k(t))+\phi_\epsilon(y_k+X_k(t)).$$

Consider a system
\begin{equation}\label{bhsdfg7}
\dot y=f^\epsilon(t,y),\quad y(0)=\hat x .\end{equation}

We have
\begin{equation}\label{4tgg}\pm f^\epsilon_k(t,x_1,\ldots,x_{k-1},\pm X_k(t),x_{k+1},\ldots,x_m)< \dot X_k(t),\quad (t,x)\in W_X.\end{equation}

Let $y^\epsilon(t)$ stand for the solution to problem (\ref{bhsdfg7}). 
Show that for all sufficiently small $\epsilon>0$ it follows that  \begin{equation}\label{evtg45}y^\epsilon(t)\ll X(t),\quad t\in I_T.\end{equation}

Indeed, fix $\epsilon\in(0,\epsilon_0)$ and assume the converse: $$\tau=\sup\{\tau ' \mid y^\epsilon(t)\ll X(t),\quad \forall t\in[0,\tau ']\}<T.$$ 
This implies that for some number $j$ and for some positive number $\delta$  we have 
\begin{equation}\label{e5gy456gy}y^\epsilon_j(\tau)=-X_j(\tau),\quad y^\epsilon_j(t)+X_j(t)<0,\quad t\in(\tau,\tau+\delta).\end{equation}
We take the sign "$-$" before $X_j$ just for definiteness,  the case $y^\epsilon_j(\tau)=X_j(\tau)$ and $y^\epsilon_j(t)>X_j(t)$ is processed in the same way. 

By formula (\ref{4tgg}) $\dot X_j(\tau)+\dot y^\epsilon_j(\tau)>0$ and the function $X_j(t)+ y^\epsilon_j(t)$ increases provided $t-\tau>0$ is small enough. This contradicts against formula (\ref{e5gy456gy}).

By the standard theorem, $\|x(\cdot)-y^\epsilon(\cdot)\|_{C(I_T)}\to 0$ as $\epsilon\to 0$. Therefore formula (\ref{evtg45}) implies the assertion of the Theorem.

The Theorem is proved.

\begin{theorem}\label{5evt5vsssdds}Suppose that $T=\infty$ and in addition to conditions of Theorem \ref{rebyh6} assume the function $f$ to be $\omega-$periodic ($\omega>0$) in $t$:
$$f(t+\omega,x)=f(t,x)$$ and $X(\omega)\ll X(0).$

Then problem (\ref{bh7}) has a solution  $\tilde x(t)\in C^1(I_T,\mathbb{R}^m)$ such that
 $\tilde x(t)\ll X(t),\quad t\in I_T$ and $\tilde x(t+\omega)=\tilde x(t)$. \end{theorem} 
\proof Ideed, this Theorem follows from the argument above. The Poincare map $\hat x\mapsto x(\omega)$
takes the convex compact set $K=\{x\in\mathbb{R}^m\mid x\ll X(0)\}$ to itself. By the Brouwer fixed point theorem there exists an initial condition $\hat x$ such that $x(\omega)=\hat x$.

\subsubsection{Non-Lipschitz Case} In this section we  assume  that $f\in C(W_X,\mathbb{R}^m)$.

The Theorem is proved.

\begin{theorem}\label{re356byh6}
Suppose that $X_k(t)>0,\quad k=1,\ldots,m,\quad t\in I_T$ and 
 for each $(t,x)\in W_X$ one has
$$\pm f_k(t,x_1,\ldots,x_{k-1},\pm X_k(t),x_{k+1},\ldots,x_m)\le \dot X_k(t),\quad\hat x\ll X(0).$$

Then problem (\ref{bh7}) has a solution  $x(t)\in C^1(I_T,\mathbb{R}^m)$ such that $x(t)\ll X(t)$. \end{theorem}
\subsubsection{Proof of Theorem \ref{re356byh6}}By the Tietze Extension Theorem, extend the functions $f_k$ to  continuous functions of an open  neighbourhood  of the set $W_T$. 

Introduce as above the function 
$$f^\epsilon_k(t,y)= f_k(t,y)-\phi_\epsilon(y_k-X_k(t))+\phi_\epsilon(y_k+X_k(t)).$$
The parameter $\epsilon$ is also chosen to fulfil equality (\ref{e5g6}). Therefore, inequality (\ref{4tgg}) is also satisfied.

Let $\{f^{n,\epsilon}(t,x)\}$ be a sequence of functions that are smooth in some open neighbourhood of the set $W_X$ and such that \begin{equation}\label{v54}\|f^{n,\epsilon}-f^\epsilon\|_{C(W_X)}\to 0\end{equation} as $n\to\infty$. For all sufficiently large $n$ these functions satisfy inequality (\ref{e5tb54}):
$$\pm f_k^{n,\epsilon}(t,x_1,\ldots,x_{k-1},\pm X_k(t),x_{k+1},\ldots,x_m)\le \dot X_k(t).$$

Thus by Theorem \ref{rebyh6} each problem
\begin{equation}\label{abbb}\dot x^{n,\epsilon}=f^{n,\epsilon}(t,x^{n,\epsilon}),\quad x^{n,\epsilon}(0)=\hat x\end{equation}has a solution $x^{n,\epsilon}(t) \in C^1(I_T),$
\begin{equation}\label{v5t5b46}\quad x^{n,\epsilon}(t)\ll X(t).\end{equation}

\begin{lemma}\label{srv4t5e}The set $U=\{ x^{n,\epsilon}(t)\}$ is relatively compact in $C(I_T)$.\end{lemma}
{\it Proof.} Indeed, by formula (\ref{v5t5b46}) the set $U$ is bounded. If we show that it is uniformly continuous then the Arzela-Ascoli Theorem implies the Lemma.

Accomplish an estimate
\begin{align}
\| x^{n,\epsilon}(t')&-x^{n,\epsilon}(t'')\|_\infty\le\Big\|
\int_{t''}^{t'}f^{n,\epsilon}(t,x^{n,\epsilon}(t))\,dt\Big\|_\infty\nonumber\\ &\le M|t'-t''|.\nonumber\end{align}
Here $\| \cdot\|_\infty$ is the standard norm in $\mathbb{R}^m$.

By formula (\ref{v54}) the constant $M$ can be chosen as follows $M=2\|f\|_{C(W_X)}.$

The Lemma is proved.

Take a subsequence  $\{x^{n_j,\epsilon_j}\}$ that is convergent to $x(t)$ in $C(I_T)$ as $n_j\to\infty,\quad \epsilon_j\to 0$. 

Passing to the same limit in the integral equation
$$x^{n_j,\epsilon_j}(t)=\hat x+\int_0^tf^{n_j,\epsilon_j}(s,x^{n_j,\epsilon_j}(s))\,ds$$
we conclude that $x(t)$ is the desired solution to problem (\ref{bh7}).

The Theorem is proved.

\begin{theorem}\label{5evt5sdds}Suppose that $T=\infty$ and in addition to conditions of Theorem \ref{re356byh6} assume the function $f$ to be $\omega-$periodic ($\omega>0$) in $t$:
$$f(t+\omega,x)=f(t,x)$$ and $X(\omega)\ll X(0).$

Then problem (\ref{bh7}) has a solution  $\tilde x(t)\in C^1(I_T,\mathbb{R}^m)$ such that
 $\tilde x(t)\ll X(t),\quad t\in I_T$ and $\tilde x(t+\omega)=\tilde x(t)$. \end{theorem} 
\proof By Theorem \ref{5evt5vsssdds} all the problems (\ref{abbb}) have $\omega-$periodic solutions $\tilde x^{n,\epsilon}(t)$ such that 
$$ \tilde x^{n,\epsilon}(t)\ll X(t).$$ By the same argument as above, the set $\{\tilde x^{n,\epsilon}(\cdot)\}$ is  relatively compact in $C(I_\omega,\mathbb{R}^m)$. Let $x_*(t)$ be an  accumulation point of this set. Then the function $$\tilde x(t)=x_*(\tau),\quad \tau\in I_\omega,\quad  t=\tau\pmod{\omega}$$ is the periodic solution.

The Theorem is proved. 

\subsection{Estimates From Below}
Now  we are again in the conditions of Section \ref{5g6rryryrty}, particularly, $f$ is a  Lipschitz function.

A notation $\mathrm{int}\, W_X$ stand for the interior of $W_X$.

\begin{theorem}\label{rbyh6}
Suppose that $X_k(t)>0,\quad k=1,\ldots,m,\quad t\in I_T$ and 
 for each $(t,x)\in W_X$ one has
\begin{equation}\label{5tb54}\pm f_k(t,x_1,\ldots,x_{k-1},\pm X_k(t),x_{k+1},\ldots,x_m)\ge \dot X_k(t).\end{equation}
And let $x(t)$ be a solution to problem (\ref{bh7}) such that $(0,\hat x)\notin \mathrm{int}\, W_X$.

Then $(t,x(t))\notin \mathrm{int}\, W_X$  for all the time of existence of the solution.
\end{theorem}

\subsubsection{Proof of Theorem \ref{rbyh6}}Our argument is almost the same as in Section \ref{5by45}. We bring it just for completeness.
 
Introduce a function 
$$g^\epsilon_k(t,y)= f_k(t,y)+\phi_\epsilon(y_k-X_k(t))-\phi_\epsilon(y_k+X_k(t)).$$
The parameter $\epsilon$ is also chosen to fulfil equality (\ref{e5g6}). Therefore, inequality (\ref{4tgg}) is also satisfied. 

Inequality (\ref{5tb54}) implies
\begin{equation}\label{5qqqtb54}\pm g^\epsilon_k(t,x_1,\ldots,x_{k-1},\pm X_k(t),x_{k+1},\ldots,x_m)> \dot X_k(t),\quad (t,x)\in W_X.\end{equation}
Let $y^\epsilon$ stand for solution to system
\begin{equation}\label{ve5}
\dot y^\epsilon=g^\epsilon(t,y^\epsilon),\quad (t_0,y^\epsilon(t_0))\notin \mathrm{int}\,W_X.\end{equation}

Let us show that $(t,y^\epsilon(t))\notin \mathrm{int}\, W_X$ for all admissible  $t>t_0$.

Assume the converse: there is a moment $\tilde t>t_0$ such that $(\tilde t,y^\epsilon(\tilde t))\in \mathrm{int}\,W_X$.
Define a parameter $\tau $ as follows
$$\tau=\inf\{\tau '\in [t_0,\tilde t]\mid (t,y^\epsilon(t))\in W_X,\quad \forall t\in[\tau ',\tilde t]\}.$$

These   imply that for some number $j$ we   have $y_j^\epsilon(\tau)=X_j(\tau)$. The case $y_j^\epsilon(\tau)=-X_j(\tau)$ is treated analogously.

Inequality (\ref{5qqqtb54}) gives $\dot y^\epsilon_j(\tau)-\dot X_j(\tau)>0$, and the function $y_j^\epsilon(t)-X_j(t)$ increases,  consequently, for some small $\delta>0$ it follows that
$$y_j^\epsilon(t)>X_j(t),\quad t\in(\tau,\tau+\delta).$$ This contradicts to the definition of $\tau$.

Now prove the Theorem. Suppose that there are  numbers $t_0< \tilde t$ such that for the solution $x(t)$ we have $(t,x(t))\in U,\quad t\in[t_0,\tilde t]$  and $$(\tilde t,x(\tilde t))\in \mathrm{int}\,W_X,\quad (t_0,x(t_0))\notin \mathrm{int}\,W_X.$$ Taking $y^\epsilon(t_0)=x(t_0)$ and approximating the solution $x(t)$ by the solutions to problem (\ref{ve5})  
$$\|x-y^\epsilon\|_{C[t_0,\tilde t]}\to 0,\quad \epsilon\to 0$$ we get the contradiction.

The Theorem is proved.

\subsection{Applications: Stability Theory}

Consider a nonlinear system
\begin{equation}\label{ev5vtesfd}
\dot x_k=\sum_{j=1}^ma_{kj}(t)x_j+\psi_k(t,x),\qquad k=1,\ldots, m. \end{equation}
Here the functions are as follows
$$a_{ij}\in C[0,\infty),\quad\psi_k\in C(A),\quad A=\{(t,x)\in\mathbb{R}^{m+1}\mid t\ge 0,\quad \|x\|_\infty\le r\},$$
and $r>0$ is a constant. The functions $\psi$ are also locally Lipschitz in the second argument and for some constants $\lambda>1,\quad c\ge 0$ it follows that
$$|\psi_k(t,x)|\le c\|x\|_\infty^\lambda,\quad (t,x)\in A.$$

Introduce a function
$$ p(t)=\max_n\Big\{a_{nn}(t)+\sum_{j\in J_n}|a_{nj}(t)|\Big\},\quad J_n=\{1,\ldots,m\}\backslash\{n\}.$$

\begin{prop}\label{ecv56}
Assume that 
$$\sup_{t\ge 0}\Big\{\int_0^tp(s)ds\Big\}<\infty\quad\mathrm{and}\quad 
c\cdot\sup_{t\ge 0}\Big\{ \int_0^te^{(\lambda-1)\int_0^\xi p(s)ds}d\xi\Big\}<\infty.$$ 
Then the zero solution to system (\ref{ev5vtesfd}) is stable in the sense of Lyapunov.\end{prop}
In the linear case ($c=0$) this proposition does not follow directly from Levinson's theorem \cite{Levinson}, but it perhaps follows from  a modification of the argument of that theorem.

To prove this Proposition \ref{ecv56} we employ Theorem \ref{rebyh6} with $X_k(t)=X(t)$. The scalar function $X$ is the solution to the following problem
$$\dot X= p(t)X+cX^\lambda,\quad X(0)=\hat X>0.$$
This Bernoulli equation is easily solved:
$$X(t)=\frac{e^{\int_0^tp(s)ds}\hat X}{\Big(1-c(\lambda-1)\hat X^{\lambda-1}\int_0^te^{(\lambda-1)\int_0^\xi p(s)ds}d\xi\Big)^{\frac{1}{\lambda-1}}}.$$
By virtue of Theorem \ref{rebyh6}, this formula implies Proposition \ref{ecv56}.

Let us introduce a function 
$$ q(t)=\min_n\Big\{a_{nn}(t)-\sum_{j\in J_n}|a_{nj}(t)|\Big\},\quad J_n=\{1,\ldots,m\}\backslash\{n\}.$$

\begin{prop}\label{ecqqv56}
Assume that 
for all $t\ge 0$ the function $q(t)$ is greater than some positive constant $C$:$$q(t)\ge C>0.$$
Then the zero solution to system (\ref{ev5vtesfd}) is  Lyapunov unstable.\end{prop}
\begin{rem}Actually this is not a solely possible conclusion from Theorem \ref{rbyh6}.
 Other one  for example is as follows. Suppose that  the expression 
$$c\cdot\limsup_{t\to\infty} \int_0^te^{(\lambda-1)\int_0^\xi q(s)ds}d\xi<\infty$$  and $$\limsup_{t\to\infty} \int_0^t q(s)ds=\infty$$ then zero solution to system (\ref{ev5vtesfd}) is also  Lyapunov unstable. \end{rem}
To prove these assertions it is sufficient to put
$X_k(t)=X(t),\quad k=1,\ldots,m$ then
$$\dot X=q(t) X-cX^\lambda,\quad X(0)=\hat X>0$$ and
$$X(t)=\frac{e^{\int_0^tq(s)ds}\hat X}{\Big(1+c(\lambda-1)\hat X^{\lambda-1}\int_0^te^{(\lambda-1)\int_0^\xi q(s)ds}d\xi\Big)^{\frac{1}{\lambda-1}}}.$$

Consider another example:
\begin{equation}\label{ev5445v}
\dot x_1=\Big(\frac{e^{2t}+1}{2}-x_1\Big)^3(1+x_2^6)+x_2^2,\quad \dot x_2=x_1^4(x_2^2-e^{2t})+\cos(x_1)x_2.\end{equation}

\begin{prop}For any initial conditions $|x_1(0)|\le 1,\quad |x_2(0)|\le 1$ system (\ref{ev5445v}) has a solution $x(t)\in C^1(I_\infty)$.
\end{prop}

Indeed, to apply Theorem \ref{rebyh6} take functions
$$X_1(t)=\frac{e^{2t}+1}{2},\quad X_2(t)=e^t$$
and note that these functions satisfy the following initial value problem
$$\dot X_1=X_2^2,\quad \dot X_2=X_2,\quad X_1(0)=X_2(0)=1.$$

\section{Appendix: Proof of Theorem \ref{e5b5j}}\label{rvbere6g}Let us note that theorem \ref{e5b5j} remains valid for the space $E$ over the field $\mathbb{C},\quad \lambda=\{\lambda_j\},\quad \lambda_j\in\mathbb{C}.$
This case is reduced to the real one by considering the realification of the space $E$ with the Schauder basis $\{e_k,ie_k\},\quad i^2=-1.$

Theorem \ref{e5b5j}  generalises the corresponding result of 
   \cite{linden} from the case of Banach spaces to the case of Fr\'echet spaces. In  the whole, our proof follows in the stream of  \cite{linden} but at several points our  argument is considerably differs from that book.
So we bring the proof for exposition's completeness sake.
\subsection{Preliminary Lemmas}

\begin{lemma}\label{sdfsdffff}
The operators $\mathcal{P}_n:E\to E$ are continuous.\end{lemma}
\proof
Consider the seminorms 
$$|x|_n=\sup_{j\in\mathbb{N}}\Big\|\sum_{k=1}^j x_ke_k\Big\|_n.$$
One evidently has $\|\cdot\|_j\le|\cdot|_j.$ So that the identity map
$$\mathrm{id}\,(E,\{|\cdot|_k\})\to (E,\{\|\cdot\|_k\})$$ is continuous.

The space $(E,\{|\cdot|_k\})$ is complete,   this fact is proved in the same manner as it is done in \cite{lust} for the case of Banach space.

By the open mapping theorem \cite{rob} the collections of seminorms $\{|\cdot|_k\}$ and $\{\|\cdot\|_k\}$ endow the space $E$ with the same topology.

Observe that the mappings $\mathcal{P}_n$ are continuous with respect to the seminorms $\{|\cdot|_k\}$.

The Lemma is proved.

\begin{lemma}\label{erterter55} Series 
 (\ref{cdvfghyu8765fghjknhbgv}) is convergent unconditionally iff one of the following conditions are fulfilled.

\begin{enumerate}
  \item\label{l1} For any sequence $\theta_k\in \{\pm 1\}$ the series  $$\sum_{k=0}^\infty \theta_kx_ke_k$$ is convergent.

  \item\label{l2} For each  $i$ and for each $\epsilon>0$ there is a natural $n$ such that for each finite set $\sigma\subset\mathbb{N},\quad \min\sigma>n$ one has
$$\Big\|\sum_{k\in\sigma }x_ke_k\Big\|_i<\epsilon.$$
  \item \label{l3} For any sequence $k_1<k_2<\ldots$ the series $$\sum_{i=1}^\infty x_{k_i}e_{k_i}$$ is convergent.
\end{enumerate}

If series  (\ref{cdvfghyu8765fghjknhbgv}) is convergent unconditionally then its sum is independent on permutation of its summands.

\end{lemma}The Banach space version of  lemma \ref{erterter55} contains in \cite{linden}.
The proof is transmitted to the case of Fr\'echet space directly.

\begin{cor}\label{sdfsddd} If  series (\ref{cdvfghyu8765fghjknhbgv}) is convergent then for any $i\in\mathbb{N}$ and for any $\epsilon>0$ there exists a number $n$ such that
$$\Big\|\sum_{k\in\sigma }x_ke_k\Big\|_i<\epsilon$$
if only $\min\sigma>n,\quad \sigma\subset\mathbb{N}$.

The set $\sigma$ is not necessarily finite.\end{cor}
Indeed, by Lemma \ref{erterter55}  choose $n_1$ such that for any finite $\sigma_1\subset\mathbb{N}$ one has 
$$\Big\|\sum_{k\in\sigma_1 }x_ke_k\Big\|_i<\epsilon/2.$$
Choose $n_2>n_1$ such that for any finite $\sigma_2\subset \mathbb{N},\quad \min\sigma_2>n_2$ it follows that
$$\Big\|\sum_{k\in\sigma_2 }x_ke_k\Big\|_i<\epsilon/4$$ and so on.
Let us put $\sigma_j=(n_j,n_{j+1}]\cap\sigma,\quad j\in\mathbb{N}.$ It remains to observe that $\sigma=\cup_j\sigma_j$ and
$$\sum_\sigma=\sum_{\sigma_1}+\sum_{\sigma_2}+\ldots.$$

Let $S=\{\pm 1\}^\mathbb{N}$ stand for the set of sequences $\theta=(\theta_1,\theta_2,\ldots),\quad \theta_k=\pm 1.$ Endow the set $S$ with the  product topology.

\begin{lemma}\label{sdfccc}The operator $\mathcal{M}_\theta:E\to E,\quad \theta\in S$ is continuous.\end{lemma}
\proof
By well-known theorem \cite{rob} it is sufficient to check that  $\mathcal{M}_\theta$ is a closed operator.

Let $x^j=\sum_{k=1}^\infty x_{jk}e_k\to x=\sum_{k=1}^\infty x_ke_k$ as $j\to\infty.$
By lemma \ref{sdfsdffff} for any $k$ it follows that $x_{jk}\to x_k$.

Suppose 
$$\mathcal{M}_\theta x^j\to z=\sum_{k=1}^\infty z_ke_k.$$
Then $\theta_kx_{jk}\to z_k$ and $z_k=\theta_kx_k$ i.e. 
$\mathcal{M}_\theta x^j\to\mathcal{M}_\theta x$.

The Lemma is proved.

\begin{lemma}\label{sdfsdf55}
For any $i\in\mathbb{N}$ there are a constant $c>0$ and a number $i'\in\mathbb{N}$ such that for all $x\in E$ the inequality holds
$$\sup_{\theta\in S}\|\mathcal{M}_\theta x\|_i\le c\|x\|_{i'}.$$
\end{lemma}
\proof
Show that for all $x\in E$ the mapping $T_x:S\to E,\quad T_x(\theta)=\mathcal{M}_\theta x $ is  continuous.
Particularly, the mapping $\theta\mapsto \|\mathcal{M}_\theta x\|_n$ is continuous.

Indeed, let $$\theta^k=\{\theta^k_j\}_{j\in\mathbb{N}}\to\theta=\{\theta_j\}_{j\in\mathbb{N}}$$
as $k\to\infty$. This implies that for any $m\in\mathbb{N}$  there is a number $K$ such that for $k>K$ one has $\theta_j^k=\theta_j,\quad j=1,\ldots,m.$

For these $k$ it follows that
$$\mathcal{M}_{\theta^k}x-\mathcal{M}_\theta x=-2\sum_{l\in\sigma_1}x_le_l+2\sum_{l\in\sigma_2}x_le_l,\quad \min\sigma_1,\,\min\sigma_2>m.$$
By Corollary \ref{sdfsddd} for all $i\in\mathbb{N}$ we have
$$\Big\|\sum_{l\in\sigma_1}x_le_l\Big\|_i,\quad \Big\|\sum_{l\in\sigma_2}x_le_l\Big\|_i\to 0$$ as $m\to \infty.$

This implies
$$\mathcal{M}_{\theta^k}x-\mathcal{M}_\theta x\to 0$$ as $k\to\infty.$

By  Tychonoff's theorem, $S$ is a compact space. Consequently, for any $i$ and $x$ it follows that
$$\sup_{\theta\in S}\|\mathcal{M}_\theta x\|_i<\infty.$$
Now the assertion of the lemma follows from the Banach-Steinhaus theorem \cite{Shvarz}.

The Lemma is proved.
\subsection{The Proof of The Theorem}
Let us show that the operator $\mathcal{M}_\lambda$ is defined for all $x\in E$. Introduce a notation
$$b_{nm}=\sum_{n\le k\le m}\lambda_kx_ke_k,\quad a_{nm}=\sum_{n\le k\le m}x_ke_k.$$

We wish to show that for each $j\in\mathbb{N}$ it follows that
$$\|b_{nm}\|_j\to 0,\quad n,m\to\infty.$$
There exists an element $f\in E^*$ such that $f(b_{ nm})=\|b_{nm}\|_j$ and $|f(x)|\le \|x\|_j,\quad x\in E$ \cite{rob}. The element $f$ depends on $n,m,j$.

Then $f(b_{nm})=\sum_{n\le k\le m}\lambda_kx_kf(e_k).$ Define a sequence $\theta\in S$ as follows. For $x_kf(e_k)\ge 0$ put $\theta_k=1$ and $\theta_k=-1$ otherwise.

Thus we have
$$\|b_{nm}\|_j\le \sup_k|\lambda_k|\sum_{n\le k\le m}\theta_kx_kf(e_k).$$
From this formula it follows that
$$\|b_{nm}\|_j \le \|\lambda\|_{\infty} f(\mathcal{M}_\theta a_{nm})\le \|\lambda\|_\infty \|\mathcal{M}_\theta a_{nm}\|_j.$$
By Lemma \ref{sdfsdf55} there is a number $i\in\mathbb{N}$ and a constant $c>0$ such that 
$$ \|\mathcal{M}_\theta a_{nm}\|_j\le c\| a_{nm}\|_i.$$ The parameters $i,c$ are independent on $a_{nm}$ and $\theta\in S$.
Since the series (\ref{cdvfghyu8765fghjknhbgv}) is convergent, $a_{nm}\to 0$ as $n,m\to\infty$ and so is $\mathcal{M}_\theta a_{nm}\to 0$. Thus $\mathcal{M}_\lambda x$ is defined for all $x\in E$ and $\lambda\in\ell_\infty.$

Now replacing $b_{nm}$ with the partial sums $b_n=\sum_{k=1}^n\lambda_kx_ke_k$ and repeating the previous argument we obtain the assertion of the theorem.

Theorem \ref{e5b5j} is proved.

\subsection*{Acknowledgement} The author wishes to thank Prof. Yu. A. Dubinskii and Prof. E. I. Kugushev for useful discussions.


\begin{thebibliography}{99}
\bibitem{Levinson} Earl A. Coddington, Norman Levinson: Theory of ordinary differential equations. New York : McGraw-Hill ; New Delhi : Tata McGraw-Hill, 1955.
\bibitem{Dubinskii}Yu. A. Dubinskii The Cauchy problem and pseudodifferential operators in the complex domain.   Russian Mathematical Surveys(1990),45(2):95.

\bibitem{kuk} S. Kuksin: Analysis of Hamiltonian PDEs.  Oxford, 2000.

\bibitem{linden}J. Lindenstrauss and L. Tzafriri. Classical Banach spaces. I. Springer-Verlag, Berlin, 1977.
\bibitem{McLeod} J. B. McLeod, On an infinite set of nonlinear differential equations, Quart. J. Math. Oxford Ser (2) 13 (19620, 119-128.

\bibitem{mill}V.M. Millionshchikov,  On the theory of differential equations in locally convex spaces, Mat. Sb. 57
(1962), 385-406. MR 27 No 6002.
\bibitem{nir} L. Nirenberg: Topics in Nonlinear Functional Analysis, New York, 1974.

  

\bibitem{rob}  A. P. Robertson and W. Robertson,
Topological Vector Spaces, Cambridge Uni-versity Press, 1964
\bibitem{sam}  A. M. Samoilenko,  Yu V. Teplinskii:
Countable Systems of Differential Equations, Brill, 2003.  

\bibitem{lust}L.A. Sobolev, V.J. Lusternik: Elements of Functional Analysis. New York, 1975.

\bibitem{Shvarz} L. Schwartz Analyse mathe'matique. Hermann, 1967. vol. 2.
\bibitem{wite} Warren H. Wite A Global Existence Theorem to Smoluchowsi's Coagulation Equation. Proceedings of the AMS Vol. 80 No 2 Oct. 1980.
\bibitem{tr_zu} D. Treschev, O. Zubelevich: Introduction to the Perturbation Theory of Hamiltonian Systems, Springer, 2003.

\bibitem{zu} O. Zubelevich On the Majorant Method for the Cauchy-Kovalevskaya Problem, Math. Notes, 69:3 (2001), 329-339.









\end{thebibliography}
 \end{document}